\documentstyle[12pt]{article}
\def\dis{\displaystyle}

\begin{document}
\hspace{.5cm}

{\sf
\centerline{\bf Integration through Generalized Sequences }

\vspace{10pt}

\centerline{Elena Toneva}

\vspace{25pt}

\centerline{Dedicated to Professor Y. Tagamlitzki}

\vspace{25pt}

{\it{Abstract:   
The process of integration was a subject of significant development during the last century.
Despite that the Lebesgue integral is complete and has many good properties, its inability to 
integrate all derivatives prompted the introduction of new approaches - Denjoy, Perron and others
introduced new ways of integration aimed at preserving the good properties of the Lebesgue integral
but extending the set of  functions to which it could be applied. The goal was achieved but neither 
of the new approaches was elegant or simple or transparent. In the 50s a new integral was introduced, 
independently by Kurzweil and Henstock,in a very 
simple, Riemann like way, but it turned out that it was more powerful than the Lebesgue integral.
There are many names attached to this integral, I will use here the name Henstock
  integral. The goal of this article is to introduce the generalized Riemann integral and the Lebesgue  integral using generalized sequences.
}}

\vspace{25pt}

Let $[a, b]$ be a compact subinterval of $R$. A function $\delta:[a,b]\rightarrow R$ is called {\it{gauge}} if it is positive.
Let $P=\{x_0,x_1,...,x_n;c_1,c_2,...,c_n\}$ be a tagged partition of $[a, b]$ ($x_0 < x_1 < ... < x_n, c_i \in [x_{i-1},x_i]$
and $x_0=a, x_n=b$). Let $f:[a,b]\rightarrow R$. With $S(P,f) = \dis{\sum_{i=1}^n f(c_i)(x_i -x_{i-1})}$ we will denote the standard
Riemann sum of $f$ corresponding to the partition $P$. A tagged partition $P$ is called $\delta$-fine if
$[x_{i-1},x_i] \subset (c_i - \delta(c_i),c_i + \delta(c_i))$ for each $i=1,2,...,n$.

\vspace{25pt}

{\bf Definition 1.} A function $f:[a,b]\rightarrow R$ is Henstock integrable on $[a, b]$ if there exists a real number $L$ with the following 
property: for each $\epsilon > 0$ there is a gauge $\delta$ such that for each $\delta$-fine partition $P$, $|S(P,f) - L| < \epsilon$.  

\vspace{25pt}

{\bf Definition 2.} A non-empty set $A$ is called directed if there is a relation $\alpha$ defined on $A$, which is reflexive, transitive and 
has the following property: if $c,d \in A$ there is $e\in A$, such that $c \alpha e$ and $d\alpha e$. If $c \alpha e$ we would say that $e$ follows
$c$. All non-empty sets of real numbers are directed by the standard inequality. 
\vspace{25pt}

{\bf Definition 3.} A function defined on a directed set is called generalized sequence.

\vspace{25pt}
The concept of generalized sequence provides a unified approach to variety of limit processes [Be]. It is well known that
the Riemann integral can be defined in terms of a generalized sequence.

Let $A=\{(P,\delta):P$ is a tagged partition of $[a,b], \delta $ is a gauge and $P$ is $\delta$-fine$\}$. The Cousin's lemma (see [B]) 
implies that $A$ is not empty. We define a relation on $A$ in the following way: we say that $ (P,\delta)$ follows $(P_1,\delta_1)$
if $\delta(x) \leq \delta_1(x)$ on the interval $[a,b]$. This relation provides a direction on $A$.

\vspace{25pt}
Let $f:[a,b]\rightarrow R$. We  define a function on $A$ in the following way:
$I(\gamma)= S(P,f)$, where $\gamma = (P,\delta)$.

\vspace{25pt}

{\bf Theorem 1.} $f$ is Henstock integrable on $[a,b]$ iff the sequence $\{I(\gamma)\}$ is convergent.

\vspace{25pt}
{\bf Proof:} (i) Assume that $f$ is Henstock integrable on $[a,b]$. We will  prove that  the 
generalized sequence is convergent and its limit is the Henstock integral of $f$ - $L$. Let $\epsilon > 0$. 
According to definition 1., there is a gauge $\delta$, such that if $P$ is $\delta$-fine, then 
$|S(P,f) - L| < \epsilon$. Let $\gamma_0 = (P_0,\delta)$, where $P_0$ is some $\delta$-fine partition.
Let $\gamma$ follows $\gamma_0$. That means that $\gamma = (Q,g)$ where $g \leq \delta$ and $Q$ is $g$-fine.
But this implies that $Q$ is also $\delta$-fine and $|S(Q,f) - L| < \epsilon$. In terms of the sequence
$\{I(\gamma)\}$ this means that we have $\gamma_0$, such that for each $\gamma$ that follows $\gamma_0$,
$|I(\gamma) - L| < \epsilon$. Since this holds for each $\epsilon$, $ \lim_{\gamma}I(\gamma)= L$.

(ii) Assume that $ \lim_{\gamma}I(\gamma)= L$. This means that for each $\epsilon > 0$ there is a $\gamma_0 \in A$,
such that if $\gamma$ follows $\gamma_0$, then $|I(\gamma) - L| < \epsilon$. Let $\gamma_0 = (P_0,\delta)$.
Let $Q$ be a $\delta$-fine partition. Then $(Q,\delta)$ will follow $\gamma_0$ and 
$|I(\gamma) - L| < \epsilon$. But $I(\gamma)= S(Q,f)$. So we found a gauge that satisfies definition 1. Q.E.D.

 \vspace{25pt}

{\bf Theorem 2. } The Henstock integral as linear and preserves the inequalities.

{\bf Proof:} Follows immediately from the properties of generalized sequences.

\vspace{25pt}

It is well known that that a function is Lebesgue integrable iff the function and its absolute value are Henstock integrable.
\vspace{25pt}

{\bf Theorem 3.}  $f$ is Lebesgue integrable on $[a,b]$ iff
both sequences $I(\gamma)$ and $J(\gamma)=S(P,|f|)$, where $\gamma=(P,\delta)$, are convergent.
\vspace{25pt}

Theorem 3. provides an opportunity to define the Lebesgue integral and Lebesgue measure in an easy and direct way.

\vspace{2cm}

{\bf References }
\vspace{25pt}

[B] Bartle, R., G., A Modern Theory of Integration, Graduate Studies in Mathematics, vol.32, MAA, 2001

\vspace{25pt}

[Be] Beardon, A., F. Limits, a New Approach to Real  Analysis,  Springer, 1997

\vspace{25pt}

[Bu] Burk, F., E., A Garden of Integrals, MAA, 2007

\vspace{25pt}

[G] Gordon, R., A., The Integrals of Lebesgue, Denjoy, Perron, and Henstock, Graduate Studies in Mathematics, vol.4, MAA, 1994

\vspace{25pt}

[H] Henstock, R., Theory of Integration,  Butterworths,  London,1963

\vspace{25pt}

[K] Kurzweil, J., Generalized Ordinary Differential Equations and Continuous Dependence on a Parameter,  Czech.  Math. J., 7 (1957), 418-446

\vspace{25pt}

[L] Leader, S., The Kurzweil-Henstock Integral and its Differentials, A Unified Theory
of  Integration on $R$ and $R^n$,  Marcel Dekker, 2001

\vspace{25pt}

[P] Petkov, P., Introduction in the Theory of Integration, Bulgarian Institute for Analysis
and Investigations, Sofia, 2000 (in Bulgarian)

\end{document}